\documentclass[authoryear, 12pt]{elsarticle}
\usepackage{fullpage}
\usepackage{graphics}
\usepackage{graphicx} 
\usepackage{amsmath}
\usepackage{amssymb}
\usepackage{url} 
\usepackage{epstopdf}

\journal{$~$}
\begin{document}
\begin{frontmatter}
\title{Spatial Process Gradients and Their Use in Sensitivity Analysis
  for Environmental Processes}
\author[ncsu]{Maria A. Terres\corref{cor1}\fnref{fn1}}
\ead{maria.terres@duke.edu}
\author[duke]{Alan E. Gelfand\fnref{fn2}}
\ead{alan@stat.duke.edu}
\cortext[cor1]{Corresponding author}
\fntext[fn1]{Maria A. Terres is Postdoctoral Research Scholar}
\fntext[fn2]{Alan E. Gelfand is Professor}
\address[ncsu]{Department of Statistics, North Carolina State University, 5109 SAS Hall, Raleigh NC 27695}
\address[duke]{Department of Statistical Science, Duke University, Box
  90251, Durham NC 27708} 

\begin{abstract} 
This paper develops methodology for local sensitivity analysis based
on directional derivatives associated with spatial processes. Formal
gradient analysis for spatial processes was elaborated in previous papers,
focusing on distribution theory for directional derivatives associated
with a response variable assumed to follow a Gaussian process
model.  In the current work, these ideas are extended to additionally
accommodate a continuous covariate whose directional derivatives are
also of interest and to relate the behavior of the directional
derivatives of the response surface to those of the covariate surface.
It is of interest to assess whether, in
some sense, the gradients of the response follow those of the
explanatory variable. The joint Gaussian structure of all variables,
including the directional derivatives, allows for explicit
distribution theory and, hence, kriging across the spatial region
using multivariate normal theory. Working within a Bayesian
hierarchical modeling framework, posterior samples enable all gradient
analysis to occur post model fitting. As a proof of concept, we show
how our methodology can be applied to a standard geostatistical
modeling setting using a simulation example. For a real data
illustration, we work with point pattern data, deferring our gradient
analysis to the intensity surface, adopting a log-Gaussian Cox process
model. In particular, we relate elevation data to point patterns associated with
several tree species in Duke Forest. 
\end{abstract}
\end{frontmatter}
\section{Introduction}
Increasingly, data is being collected at geo-referenced locations.
For a region of interest $D$, the set of conceptual responses
$\{Y({\bf s}):{\bf s}\in D\}$ can be viewed as a realization of a
random surface, observed at a finite set of locations. While covariate
information may explain a substantial portion of the variation in
response, there is often underlying spatial structure that is
difficult to measure. Inference on this spatial structure can be made
via the parameters in a spatial process model. Under these models,
prediction at unobserved locations, or kriging, is possible, enabling
interpolation across the region. 

Spatial regression models commonly assume a linear relationship and
make inference based on the coefficient assigned to the
covariate. This coefficient describes the expected change in response
given a unit change in covariate, thus providing a global measure for the
sensitivity of the response to the covariate.  However, it is expected
that the relationship between the variables may vary locally over the study
region.  Such local, or second-order, behavior
can be studied through spatial sensitivity or gradient analysis.  

A spatial gradient analysis will enable spatial examination of a response
variable's sensitivity to a covariate across the region of
interest. The sensitivity of the response variable may vary based on
the rate of change for the covariate or due to additional unaccounted
for factors, resulting in areas where the relationship appears weaker
or stronger. 
Models allowing for spatially varying coefficients, such as \cite[]{
  fotheringham2002geographically, Gelfand2003}, provide some 
similar inference but assume a more complex model structure. The
methodology proposed here assumes a standard spatial linear regression
model but provides a post model fitting framework for
examining the variation in the response's sensitivity to the covariate.

In ecology such sensitivities are typically
discussed when relating plant characteristics to climate. For instance,
researchers are increasingly interested in characterizing how abundance
and frequency of tree species relate to changes in variables such as
temperature and precipitation, in order to learn about the expected effects of
climate change on range distributions \cite[e.g.][]{Thuiller2004, CanhamThomas2010, Thomas2010}. 
These analyses focus on a comparison of sensitivities across species, 
while little
consideration is given to how sensitivities may vary spatially within
a given species. Our approach is aimed at the latter question. 

The interpretation of the coefficient in a spatial regression as a
global gradient, $dE(Y({\bf s}))/$\linebreak$dX({\bf s})$, inspires
consideration of local sensitivities through directional
derivative processes, or spatial gradients. Spatial gradients under Gaussian
processes were elaborated by \cite{Banerjee2003} to address the rate
of change of a spatial surface at a given point in a given
direction. Their paper defines directional derivative processes with
corresponding distribution theory to enable interpolation across a
region. The gradient distributions are fully determined by the spatial
model parameters, allowing all gradient analysis to occur post model
fitting. Distributions for derivatives of Gaussian processes have also
been discussed  in the context of observed derivatives of functions
\cite[e.g.][]{o1992some, murray2003derivative}, as well as for random
fields more generally \cite[]{adler1981geometry}.
In all of the previous work with Gaussian process directional derivatives
\cite[e.g.][]{BanerjeeGelfand2003, BanerjeeGelfand2006, Majumdar2006},
the researchers have considered the rates of change of a response
surface with the mean surface modeled as a 
linear function of a set of fixed covariates. In contrast, to accommodate the
desired spatial sensitivity analysis, we will assume a single covariate of
interest whose surface is spatially smooth such that it too can be
treated as a realization of a stochastic process. The behavior of the
response and covariate processes, as well as their associated
derivative processes, are then considered jointly and functions of
these derivatives can be explored. 

The contribution of this paper is to extend the existing spatial gradient
theory to accommodate spatial sensitivity analysis by modeling the response
and covariate jointly. Working within a hierarchical Bayesian modeling
framework, corresponding gradients for the spatial
surfaces can be sampled simultaneously from the joint predictive
distribution post model fitting. Under a significant regression
relationship it is not sensible to investigate the gradient behavior
of the surfaces marginally. Suitable comparison between the gradient surfaces
illustrates how sensitive the response surface is to the covariate
surface, as well as the strength of this relationship. The former is
accomplished through comparison between the directions of the maximum
gradient at a given location; the latter requires consideration of
their directional derivatives relative to one another. In particular,
we introduce two new spatial processes, a local directional
sensitivity process and a spatial angular discrepancy process. These
inferential tools are developed and carried out on simulated data in
the context of a customary geostatistical model \cite[]{BCG2003,
  cressie2011statistics} as 
well as with an ecological dataset where we connect point patterns of
trees with elevation. 

In Section \ref{sec:dists} the formal distribution theory for the
spatial gradients is extended to the multivariate case. Section
\ref{sec:model} outlines the modeling framework for our examples. This
section also defines the two processes of interest, namely the local
directional sensitivity process and the spatial angular discrepancy
process. Section \ref{sec:sim} provides a simulated example with a
multivariate Gaussian process setup as a proof of concept. Section
\ref{sec:data} provides an analysis of point pattern data from Duke
Forest, extending the analysis techniques to a non-Gaussian response;
the intensity of the point pattern, modeled using a log-Gaussian Cox
process, is explored by employing a spatial gradient chain
rule. Finally, Section \ref{sec:summary} summarizes the contributions
of the paper and suggests future work. 

\section{Distribution Development}
\label{sec:dists} In this section we review the definitions and distributions presented in \cite{Banerjee2003} and extend these
ideas to consider a multivariate Gaussian process. We assume locations ${\bf s} \in \mathbb{R}^2$, 2-dimensional Euclidean space,
however extension to a generic $d$-dimensional setting is straightforward. The process is assumed, for convenience, to be
(weakly) stationary such that the covariance function, $Cov(Y({\bf s}), Y({\bf s}'))$, depends only on the separation vector
${\boldsymbol\delta}={\bf s}-{\bf s}'$. In fact, in our examples we adopt isotropic covariance functions that depend only on the
length of the separation vector, $||{\boldsymbol\delta}||$.

Consider two surfaces $\{(Y({\bf s}),X({\bf s})): {\bf s}\in \mathbb{R}^2\}$ drawn from a joint Gaussian process specified such
that $X({\bf s})$ has constant mean, say $\alpha_0$, and covariance function $G({\boldsymbol\delta})$. Given $X({\bf s})$,
$Y({\bf s})$ has mean $\beta X({\bf s})$ and covariance function $K({\boldsymbol\delta})$. Observed at a set of locations ${\bf
Y}=(Y({\bf s}_1), \hdots ,Y({\bf s}_n))$, we write: ${\bf Y}|{\bf X}  \sim N(\beta {\bf X}, K(\cdot))$ with
  ${\bf X} \sim N(\alpha_0, G(\cdot))$. Considered jointly, we have:
\begin{align*}
\left(\begin{array}{c}
{\bf Y}\\ {\bf X}\end{array} \right)
&\sim N\left( \left(\begin{array}{c}\alpha_0\beta{\bf 1}\\ \alpha_0{\bf 1}\end{array}\right),
\left( \begin{array}{cc}K(\cdot) + \beta^2G(\cdot) & \beta G(\cdot)\\
    \beta G(\cdot) & G(\cdot)\end{array}\right)
\right)
\end{align*}
where $G(\cdot)$ and $K(\cdot)$ are matrices of the covariance
functions with entry $i,j$ evaluated at ${\boldsymbol\delta}={\bf
  s}_i-{\bf s}_j$.

We follow the notation and theory in \cite{Banerjee2003}. Suppose mean square differentiable processes $Y({\bf s})$ and $X({\bf s})$. That
is, for $Y(\mathbf{s})$, at $\mathbf{s}_{0}$ there exists a vector $\nabla_{Y}(\mathbf{s}_{0})$ such that for any scalar $h$ and
any unit vector $\mathbf{u}$, $Y(\mathbf{s}_{0} + h\mathbf{u}) = Y(\mathbf{s}_{0}) + h \mathbf{u}^{T}\nabla_{Y}(\mathbf{s}_{0})+
r(\mathbf{s}_{0},h\mathbf{u})$ where $r(\mathbf{s}_{0},h\mathbf{u}) \rightarrow 0$ in the $L_{2}$ sense as $h \rightarrow 0$.
Similarly, for $X(\mathbf{s})$.

In particular, define the finite difference processes at scale $h$ in direction ${\bf u}$:
\begin{align*}
Y_{{\bf u},h}({\bf s})=\frac{Y({\bf s}+h{\bf u})-Y({\bf s})}{h}\\
X_{{\bf u},h}({\bf s})=\frac{X({\bf s}+h{\bf u})-X({\bf s})}{h}
\end{align*}
where ${\bf u}$ is a unit vector. Taking the limit as $h$ tends to 0, \cite{Banerjee2003} define the directional derivative
processes in the direction ${\bf u}$:
\begin{align*}
D_{\bf u}Y({\bf s}) &= \lim_{h\to 0}Y_{{\bf u},h}({\bf s})= {\bf
  u}'\nabla_Y({\bf s})\\ 
D_{\bf u}X({\bf s}) &= \lim_{h\to 0}X_{{\bf u},h}({\bf s})= {\bf u}'\nabla_X({\bf s})
\end{align*}
where $\nabla_X({\bf s})=(D_{{\bf e}_1}X({\bf s}), D_{{\bf e}_2}X({\bf s}))'$ is the vector of directional derivatives in the
orthonormal basis directions ${\bf e}_1=(1,0)$ and ${\bf e}_2=(0,1)$ for $\mathbb{R}^2$. We can study the directional derivative
processes for any ${\bf u}$ by working with the basis set $\nabla_Y({\bf s})$ and $\nabla_X({\bf s})$.

From \cite{Banerjee2003}, we know that if $Y({\bf s})$ and $X({\bf s})$ are
stationary Gaussian processes, then the resulting marginal
distributions involving the directional derivatives will be stationary
Gaussian processes as well. Note, isotropy in the $Y({\bf s})$ process does
not induce isotropy in the $D_{\bf
  u}Y({\bf s})$ process; only stationarity will be inherited.  Similar to the
discussion in their paper, we know by
linearity that $(Y({\bf s}),X({\bf s}),Y_{{\bf
    u},h}({\bf s}),X_{{\bf u},h}({\bf s}))'$ will be a stationary multivariate Gaussian
process. And then, by a standard limiting moment generating function argument, $(Y({\bf s}),X({\bf s}),D_{\bf u}Y({\bf s}),D_{\bf u}X({\bf s}))'$ will
also be a stationary multivariate Gaussian process.

To explicitly provide the joint distribution, we derive the cross covariance structure by examining pair-wise covariances between
the response and covariate processes and their directional derivatives.  For notational convenience, write the marginal
covariance function of $Y({\bf s})$ to be $\tilde{K}(\cdot)=K(\cdot) + \beta^2G(\cdot)$. Assume the $Y({\bf s})$ and $X({\bf s})$ processes are mean zero,
setting $\alpha_0=0$, since in practice the gradients are calculated for the mean zero residual process. If $E(Y({\bf s}))=0$, then
$E(D_{\bf u}Y({\bf s})) = 0$, so the joint processes will all be mean zero. We calculate the covariances associated with the
directional derivatives by taking the limits of the covariances
corresponding to the analogous finite difference process.

The covariances for the response surface are derived in \cite{Banerjee2003}:
\begin{align*}
Cov(Y_{{\bf u},h}({\bf s}),Y_{{\bf u},h}({\bf s}'))
& = \frac{2\tilde{K}({\boldsymbol\delta})-\tilde{K}({\boldsymbol\delta}+h{\bf u})-\tilde{K}({\boldsymbol\delta}-h{\bf u})}{h^2}\\
Cov(D_{\bf u}Y({\bf s}), D_{\bf u}Y({\bf s}')) & = \lim_{h\to0}Cov(Y_{{\bf u},h}({\bf s}),Y_{{\bf u},h}({\bf s}')) =-{\bf u}'\Omega_{\tilde{K}}{\bf u}\\
\end{align*}
\begin{align*}
Cov(Y({\bf s}),Y_{{\bf u},h}({\bf s}')) & = \frac{\tilde{K}({\boldsymbol\delta}-h{\bf u})-\tilde{K}({\boldsymbol\delta})}{h}\\
Cov(Y({\bf s}), D_{\bf u}Y({\bf s}')) & = \lim_{h\to0}Cov(Y({\bf s}),Y_{{\bf u},h}({\bf s}')) = D_{\bf u}\tilde{K}(-{\boldsymbol\delta})\\
\end{align*}
and covariances for the covariate surface are analogous:
\begin{align*}
Cov(X_{{\bf u},h}({\bf s}),X_{{\bf u},h}({\bf s}')) & = \frac{2G({\boldsymbol\delta})-G({\boldsymbol\delta}+h{\bf u})-G({\boldsymbol\delta}-h{\bf u})}{h^2}\\
Cov(D_{\bf u}X({\bf s}), D_{\bf u}X({\bf s}')) & = \lim_{h\to0}Cov(X_{{\bf u},h}({\bf s}),X_{{\bf u},h}({\bf s}')) = -{\bf u}'\Omega_{G}{\bf u}\\
\end{align*}
\begin{align*}
Cov(X({\bf s}),X_{{\bf u},h}({\bf s}')) & = \frac{G({\boldsymbol\delta}-h{\bf u})-G({\boldsymbol\delta})}{h}\\
Cov(X({\bf s}), D_{\bf u}X({\bf s}')) & = Cov(X({\bf s}),X_{{\bf u},h}({\bf s}')) = D_{\bf u}G(-{\boldsymbol\delta})\\
\end{align*}
To fully describe the joint distribution we derive the covariances between
response and covariate surfaces similarly:
\begin{align*}
Cov(Y({\bf s}),X_{{\bf u},h}({\bf s}')) & = \frac{\beta G({\boldsymbol\delta}-h{\bf u})-\beta G({\boldsymbol\delta})}{h}\\
Cov(Y({\bf s}), D_{\bf u}X({\bf s}')) & = \lim_{h\to0}Cov(Y({\bf s}),X_{{\bf u},h}({\bf s}')) = \beta D_{\bf u}G(-{\boldsymbol\delta})\\
\end{align*}
\begin{align*}
Cov(X({\bf s}),Y_{{\bf u},h}({\bf s}')) & = \frac{\beta G({\boldsymbol\delta}-h{\bf u})-\beta G({\boldsymbol\delta})}{h}\\
Cov(X({\bf s}), D_{\bf u}Y({\bf s}')) & = \lim_{h\to0}Cov(X({\bf s}),Y_{{\bf u},h}({\bf s}')) = \beta D_{\bf u}G(-{\boldsymbol\delta})\\
\end{align*}
\begin{align*}
Cov(X_{{\bf u},h}({\bf s}), Y_{{\bf u},h}({\bf s}') & = \frac{2\beta G({\boldsymbol\delta})-\beta G({\boldsymbol\delta}+h{\bf u}) -\beta G({\boldsymbol\delta}-h{\bf u})}{h^2}\\
Cov(D_{\bf u}X({\bf s}), D_{\bf u}Y({\bf s}')) & = \lim_{h\to0}Cov(X_{{\bf u},h}({\bf s}), Y_{{\bf u},h}({\bf s}') = -\beta {\bf u}'\Omega_G({\boldsymbol  \delta}) {\bf u}\\
\end{align*}
where $(\Omega_G({\boldsymbol\delta}))_{ij}=\partial ^2G({\boldsymbol\delta})/\partial
{\boldsymbol\delta}_i\partial{\boldsymbol\delta}_j$ and $D_{\bf u}G({\boldsymbol
  \delta})=\lim_{h\to0}(G({\boldsymbol\delta} -h{\bf u})-G({\boldsymbol
  \delta}))/h$.

Relationships between the response surface, the covariate surface and
their corresponding directional derivative surfaces can be described
through the 6-dimensional multivariate stationary Gaussian process
${\bf Z}({\bf s})=(Y({\bf s}), X({\bf s}), \nabla_Y({\bf s}),
\nabla_X({\bf s}))'$.
Using the covariances calculated above, the associated
cross-covariance matrix for ${\bf Z}$ will be:
\begin{align*}
V_{\bf Z}({\boldsymbol\delta}) & =
\left( \begin{array}{cccc}\tilde{K}({\boldsymbol\delta}) & \beta G({\boldsymbol\delta}) & -\nabla\tilde{K}({\boldsymbol\delta})' & -\beta \nabla G({\boldsymbol\delta})'\\
    \beta G({\boldsymbol\delta}) & G({\boldsymbol\delta}) & -\beta \nabla G({\boldsymbol\delta})' & - \nabla G({\boldsymbol\delta})'\\
    \nabla\tilde{K}({\boldsymbol\delta}) & \beta\nabla G({\boldsymbol\delta}) & -H_{\tilde{K}}({\boldsymbol\delta}) & -\beta H_G({\boldsymbol\delta})\\
    \beta \nabla G({\boldsymbol\delta}) & \nabla G({\boldsymbol\delta}) & -\beta H_G({\boldsymbol\delta}) & -H_G({\boldsymbol\delta})
    \end{array}\right)\\~&~\\&=
\left( \begin{array}{cccc}K({\boldsymbol\delta}) + \beta^2G({\boldsymbol\delta}) & \beta G({\boldsymbol\delta}) & -\nabla K({\boldsymbol\delta})' -\beta^2\nabla G({\boldsymbol\delta})' & -\beta \nabla G({\boldsymbol\delta})'\\
    \beta G({\boldsymbol\delta}) & G({\boldsymbol\delta}) & -\beta \nabla G({\boldsymbol\delta})' & - \nabla G({\boldsymbol\delta})'\\
    \nabla K({\boldsymbol\delta}) +\beta^2\nabla G({\boldsymbol\delta}) & \beta\nabla G({\boldsymbol\delta}) & -H_K({\boldsymbol\delta})-\beta^2H_G({\boldsymbol\delta}) & -\beta H_G({\boldsymbol\delta})\\
    \beta \nabla G({\boldsymbol\delta}) & \nabla G({\boldsymbol\delta}) & -\beta H_G({\boldsymbol\delta}) & -H_G({\boldsymbol\delta})
    \end{array}\right)
    \end{align*}
where $\nabla K({\boldsymbol\delta})$ is a $2\times1$ gradient vector associated
with $K({\boldsymbol\delta})$, and $H_K({\boldsymbol\delta})$ is the $2\times2$ Hessian matrix
associated with $K({\boldsymbol\delta})$.

For ${\boldsymbol\delta}={\bf 0}$, we have a block diagonal local covariance matrix:
\begin{align*}
V_{\bf Z}({\bf 0}) &=
\left( \begin{array}{cccc}K({\bf 0}) + \beta^2G({\bf 0}) & \beta G({\bf 0}) & {\bf 0}' & {\bf 0}' \\
    \beta G({\bf 0}) & G({\bf 0}) & {\bf 0}' & {\bf 0}'\\
    {\bf 0} & {\bf 0} & -H_K({\bf 0})-\beta^2H_G({\bf 0}) & -\beta H_G({\bf 0})\\
    {\bf 0} & {\bf 0} & -\beta H_G({\bf 0}) & -H_G({\bf 0})
    \end{array}\right).
    \end{align*}
Thus, at a location ${\bf s}$, the directional derivative surfaces will be correlated with one another, but neither will be
correlated with either of the data surfaces. Intuitively, this makes sense since we would not expect the level of the surface at
a location to be correlated with the rate of change at that location. Of course, since $(X(\mathbf{s}), Y(\mathbf{s}))'$ is a
bivariate Gaussian process, the correlation between the rates of changes is not surprising.

The Mat\'{e}rn covariance is adopted below. It depends on a smoothness parameter $\nu$ which directly controls the mean square
differentiability of process realizations \cite[]{Stein1999}. This is convenient since, again, the $Y({\bf s})$ and $X({\bf s})$
processes must be mean square differentiable for their associated directional derivative processes to be well defined.
 If we let $K(\cdot)$ and $G(\cdot)$ be Mat\'{e}rn with $\nu >1$
 then they are once (but not twice) mean square differentiable, and, if $\nu =3/2$,
the covariance functions are of the closed form $\sigma^2(1+\phi||{\boldsymbol\delta}||)\exp(-\phi||{\boldsymbol\delta}||)$. We
denote the parameters of $K(\cdot)$ as $\sigma^2_y$ and $\phi_y$, and the parameters of $G(\cdot)$ as $\sigma^2_x$ and $\phi_x$.

Under the Mat\'{e}rn covariance the components of the cross-covariance matrix will be $\nabla K({\boldsymbol \delta}) =
-\sigma_y^2\phi_y^2\exp(-\phi_y||{\boldsymbol\delta}||){\boldsymbol\delta}$, $(H_K({\boldsymbol\delta}))_{ii}=-\sigma_y^2\phi_y^2
\exp(-\phi_y||{\boldsymbol\delta}||)(1-\phi_y\delta_i^2/||{\boldsymbol\delta}||)$,
$(H_K({\boldsymbol\delta}))_{ij}$\linebreak$
=\sigma_y^2\phi_y^3
\exp(-\phi_y||{\boldsymbol\delta}||)\delta_i\delta_j/||{\boldsymbol\delta}||$,
and similar for $G(\cdot)$. Then, we have for ${\boldsymbol\delta}={\bf 0}$:
\begin{align}\label{eq:cov}
V_{\bf Z}({\bf 0}) &=
\left( \begin{array}{cccc}\sigma^2_y + \beta^2\sigma^2_x & \beta \sigma^2_x & {\bf 0}' & {\bf 0}' \\
    \beta \sigma^2_x & \sigma^2_x & {\bf 0}' & {\bf 0}'\\
    {\bf 0} & {\bf 0} & (\sigma^2_y\phi^2_y+\beta^2\sigma^2_x\phi^2_x)I_2 & (\beta \sigma^2_x\phi^2_x)I_2\\
    {\bf 0} & {\bf 0} & (\beta \sigma^2_x\phi^2_x)I_2 &(\sigma^2_x\phi^2_x)I_2
    \end{array}\right)
\end{align}
where $I_2$ is the $2\times 2$ identity matrix. As above, $D_{{\bf e}_1}Y({\bf s})$ and $D_{{\bf e}_1}X({\bf s})$ will be
correlated with one another, and similarly $D_{{\bf e}_2}Y({\bf s})$ and $D_{{\bf e}_2}X({\bf s})$ will be correlated with one
another, but all other pairings of the directional derivatives will be uncorrelated.

\section{Model Fitting and Inference}
\label{sec:model}
\subsection{Sampling Method}
Following the modeling of the previous section, let $K({\boldsymbol\delta})=\sigma_y^2\rho_y({\boldsymbol\delta})$ and
$G({\boldsymbol\delta})=\sigma^2_x\rho_x({\boldsymbol\delta})$, where the $\rho_x$  and $\rho_y$ are valid
two-dimensional correlation functions. We work with the Mat\'{e}rn class of covariance functions parameterized by $\phi$ and
$\nu$ with $\nu >1$.

 Let $ {\boldsymbol\theta}=(\alpha_0,
\beta_0, \beta_1, \sigma_x^2, \sigma_y^2, \phi_x, \phi_y, \nu_x, \nu_y)$. For locations ${\bf s_1},\hdots,{\bf s_n}$, the overall
likelihood can be written in terms of the conditional likelihoods
\begin{align*}
L({\bf {\boldsymbol\theta}; Y, X}) &\propto L({\bf Y} |
\boldsymbol\theta, {\bf X}) L({\bf X}|\boldsymbol\theta)\\
& \propto (\sigma_x^2\sigma_y^2)^{-n/2}|R_x(\phi_x,\nu_x)|^{-1/2}|R_y(\phi_y,\nu_y)|^{-1/2}\\
& \times \exp\left\{ -\frac{1}{2\sigma_x^2}({\bf X} - \alpha_0{\bf
    1})'R_x^{-1}(\phi_x,\nu_x)({\bf X} - \alpha_0{\bf 1})\right\}\\
& \times \exp\left\{ -\frac{1}{2\sigma_y^2}({\bf Y} - (\beta_0{\bf
    1} + \beta_1{\bf X}))'R_y^{-1}(\phi_y,\nu_y)({\bf Y} - (\beta_0{\bf
    1} + \beta_1{\bf X}))\right\}
\end{align*}
where ${\bf Y}=(Y({\bf s_1}),\hdots, Y({\bf s_n}))'$, $(R_{x}(\phi_x,\nu_x))_{ij}=\rho_{x}({\bf s}_i-{\bf s}_j;\phi_x,\nu_x)$ and $(R_{y}(\phi_y,\nu_y))_{ij}=\rho_{y}({\bf s}_i-{\bf
s}_j;\phi_y,\nu_y)$. The likelihood could be equivalently written in its joint form, but the conditional form is
more conducive to
interpreting and implementing the gradient analysis.

We see that we have a low dimensional parametric model, with ${\boldsymbol\theta}$ only $9$ dimensional.  We utilize fairly
non-informative priors for its components. For example, vague normal priors on $(\alpha_0, \beta_0, \beta_1)$, vague inverse
Gamma priors on $(\sigma_x^2, \sigma_y^2)$, vague Gamma priors on $(\phi_x, \phi_y)$, and $U(1,2)$ priors on $(\nu_x,\nu_y)$. The
prior on $\nu$ follows the suggestion of \cite{Stein1999} and others who observe that distinguishing $\nu=2$ from $\nu>2$ would
be very difficult in practice. This model is straight forward to fit in its conditional form, for example using the
\texttt{`spBayes'} package in \texttt{R} \cite[]{Finley2010}. Thus, assume we now have posterior samples ${\boldsymbol\theta}_l^*,
l=1,\hdots,L$, from $f({\boldsymbol\theta}|{\bf Y,X})$.

Once we  have posterior samples of the parameters, we draw samples of the gradient vectors using composition since the posterior
predictive distribution $f(\nabla_Y,\nabla_X|{\bf Y,X})=\int f(\nabla_Y,\nabla_X|{\bf
  Y,X},\boldsymbol\theta)f(\boldsymbol\theta|{\bf Y,
  X})d\boldsymbol\theta$. The cross-covariance matrix derived earlier
allows us to immediately write the joint multivariate normal
distribution given $\boldsymbol\theta$, which can be evaluated at each sample
$\boldsymbol\theta_l^*$. Based on this joint distribution,
standard multivariate normal theory allows us to write down the
desired conditional distributions needed to draw from the predictive
distribution.

For an unobserved location ${\bf s}_0$, obtaining draws of the gradient vectors is again done via the predictive distribution.
The cross covariance matrix derived enables us to write the joint distribution, from which we can derive the conditional distribution
$f(\nabla_Y({\bf s}_0), \nabla_X({\bf s}_0) | {\bf Y,X}, \boldsymbol\theta)$. If interest is also in the values of the $Y({\bf
s})$ and $X({\bf s})$ surfaces at the new location, we would
 derive the conditional distribution $f(Y({\bf s}_0), X({\bf s}_0),
 \nabla_Y({\bf s}_0),
\nabla_X({\bf s}_0) | {\bf Y,X}, \boldsymbol\theta)$, which is again straight forward given the cross covariance matrix and
allows joint prediction of the surfaces and their gradients at the new location.

If we want $Y({\bf s})$ and $X({\bf s})$ to be adjusted based on some fixed covariates, then we simply introduce such covariates
into the mean functions of the model. We create a spatial random effects model:
\begin{align*}
 Y({\bf s})|X({\bf s}) & = \beta_0 + \beta_1 X({\bf s}) +{\bf T}_y({\bf s})'\boldsymbol\gamma_y +
w_y({\bf s}) + \epsilon({\bf s})\\
X({\bf s}) & = \alpha_0 + {\bf T}_x({\bf s})'\boldsymbol\gamma_x + w_x({\bf s})
\end{align*}
where ${\bf T}_y({\bf s})$ and ${\bf T}_x({\bf s})$ are vectors of covariates used to explain the $Y({\bf s})$ and $X({\bf s})$
surfaces respectively, with coefficients $\boldsymbol\gamma_y$ and $\boldsymbol\gamma_x$; $w_y({\bf s})$ and $w_x({\bf s})$ are
independent mean-zero stationary Gaussian processes with parameters $\sigma_x^2, \sigma_y^2,\phi_x, \phi_y,\nu_x, \nu_y$ as
before; and $\epsilon({\bf s})$ is a Gaussian white-noise process with variance $\tau^2$ intended to capture measurement error or
microscale variability in the response. The $
  X({\bf s})$ process
is assumed to be a fully spatial model (no nugget effect), such as
might be used for elevation, temperature, or pollutant level.

With $ {\boldsymbol\theta}$ now extended to $\boldsymbol{\theta}=(\alpha_0, \beta_0, \beta_1, \boldsymbol\gamma_x,
\boldsymbol\gamma_y, \sigma_x^2, \sigma_y^2, \phi_x, \phi_y, \nu_x, \nu_y)$, the likelihood for locations ${\bf s_1},\hdots,{\bf
s_n}$ above can be trivially revised.

Prior selection for the parameters will be similar to the previous example. Again, this model can be implemented using
\texttt{`spBayes'}, and draws of the gradients will rely on the posterior predictive distribution, which can be calculated as
before using the derived cross-covariance matrix.

\subsection{Local Directional Sensitivity Process}
\label{subsec:direcderiv}
At a given location ${\bf s}$ there may additionally be interest in the
ratio of directional derivatives, \linebreak$D_{\bf u}Y({\bf s})/D_{\bf
    u}X({\bf s})$, corresponding to the relative rates of change in the two
surfaces in direction ${\bf u}$. This quantity is analogous to $dy/dx$ in more standard calculus applications as well as
sensitivity functions studied in sensitivity analysis \cite[]{TomovicVukobratovic1972}. With this in mind, we refer to the
resulting spatial process as the {\it local directional sensitivity process}. The choice of ${\bf u}$ will depend on the
application being considered.  For example, this direction may correspond to latitudinal direction, an elevation direction, or to
an environmental feature expected to impact the response. Large values of this process would suggest that the change in covariate
surface has a high impact on the the response surface.

The joint process defined in Section \ref{sec:dists} can be equivalently represented as a spatial random effects model (excluding
nugget effects):
\begin{align*}
Y({\bf s})|X({\bf s}) & = \beta_0 + \beta_1X({\bf s}) + w_y({\bf s})\\
X({\bf s}) & = \alpha_0 + w_x({\bf s})\\
w_y({\bf s}) &\sim GP(0,K({\boldsymbol\delta}))\\
w_x({\bf s}) &\sim GP(0,G({\boldsymbol\delta}))
\end{align*}
where $w_y({\bf s})$ and $w_x({\bf s})$ are independent processes. With this notation we can write the unconditional response surface and
corresponding directional derivative process as follows:
\begin{align*}
Y({\bf s}) & = \beta_0 + \beta_1\alpha_0 + \beta_1w_x({\bf s}) + w_y({\bf s})\\
D_{\bf u}Y({\bf s}) & = \beta_1 D_{\bf u}w_x({\bf s}) + D_{\bf u}w_y({\bf s})
\end{align*}
The local directional sensitivity process can then be written
\begin{align}\label{LDSproc}
\frac{D_{\bf u}Y({\bf s})}{D_{\bf u}X({\bf s})}&=\frac{\beta_1 D_{\bf u}w_x({\bf s}) + D_{\bf u}w_y({\bf s})}{D_{\bf u}w_x({\bf s})} = \beta_1 + \frac{D_{\bf u}w_y({\bf s})}{D_{\bf u}w_x({\bf s})}
\end{align}
We see that the multiplicative parameter, $\beta_1$, defining the overall relationship between the $X({\bf s})$ and $Y({\bf s})$
processes serves to center the local directional sensitivity
process.

As mentioned in Majumdar et al. 2006, if we consider $Y({\bf s})=\beta_0+\beta_1X({\bf s})+\epsilon({\bf s})$, then one could
write $\beta_1=dE(Y({\bf s}))/dX({\bf s})$; i.e., $\beta_1$ is describing the rate of change in $E(Y({\bf s}))$ relative to
changes in $X({\bf s})$. Again, at any location the directional derivative ratio will be centered at the global (non-directional)
derivative ratio $dE(Y({\bf
  s}))/dX({\bf s})$ plus
some (directional) spatial noise.
In this way, the directional derivative ratio process is describing the
 spatial variation in the relative rates of change between $X({\bf s})$
and $Y({\bf s})$. This is analogous to modeling adopting spatially varying coefficients, $\beta({\bf s})$, \cite[]{Gelfand2003}
but is arguably a simpler  context since the derivatives require no additional model fitting. In addition, consideration of
directional perspectives using this gradient approach allows for inference distinct from what one can learn from a
non-directional $\beta({\bf s})$ parameter.

By noting that $D_{\bf u}w_y({\bf s})/D_{\bf u}w_x({\bf s})$ is a ratio of independent mean zero normal random variables, at each
location the directional spatial noise is a Cauchy random variable with scale equal to the ratio of the respective standard
deviations, $SD(D_{\bf u}w_y({\bf s}))/SD(D_{\bf u}w_x({\bf s}))$. If ${\bf u}=(u_1, u_2)$, then $SD(D_{\bf
  u}w_y({\bf s}))=\sqrt{u_1^2(-H_K(0)_{1,1})+u_2^2(-H_K(0)_{2,2})}$. If $Y({\bf s})$ is
isotropic, then $H_K(0)=c_oI_2$ \cite[]{Banerjee2003} and
$SD(D_{\bf
  u}w_y({\bf s})) = \sqrt{u_1^2+u_2^2}\sqrt{-c_0}=\sqrt{-c_0}$. If $K(\cdot)$
and $G(\cdot)$ are Mat\'{e}rn with $\nu=3/2$, then $SD(D_{\bf u}w_i({\bf s})) =\sigma_i\phi_i$, and the scale for the Cauchy
distribution will be $\sigma_y\phi_y/\sigma_x\phi_x$.  When the respective standard deviations are equal, the scale will be 1 and
the directional derivative ratio will have a standard Cauchy distribution.

In fact, the collection of directional derivative ratios form a well defined spatial stochastic process which would naturally be
called a spatial Cauchy process (see Appendix A for details). For any set of locations ${\bf s}_1, \hdots, {\bf s}_n$ the joint
distribution is well defined. For example, consider two locations ${\bf s}$ and ${\bf s}'$, simplifying the notation for clarity:
\begin{align*}
P(\frac{D_{\bf u}w_y({\bf s})}{D_{\bf u}w_x({\bf s})}<r_1, \frac{D_{\bf
    u}w_y({\bf s}')}{D_{\bf u}w_x({\bf s}')}<r_2) & = P(\frac{n_1}{m_1}<r_1,~
\frac{n_2}{m_2}<r_2)\\
& = P(n_1<r_1m_1,~ n_2<r_2m_2,~ m_1>0,~ m_2>0)\\
& +P(n_1<r_1m_1,~ n_2>r_2m_2,~ m_1>0,~ m_2<0)\\
& + P(n_1>r_1m_1,~ n_2<r_2m_2,~ m_1<0,~ m_2>0)\\
&+ P(n_1>r_1m_1,~ n_2>r_2m_2,~ m_1<0,~ m_2<0)
\end{align*}
In turn, each of these terms can be computed as an integral involving
normal densities. For example, the first term can be written as
follows:
\begin{align}\label{eq:cauch}
 P(n_1<r_1m_1, n_2<r_2m_2, m_1>0, m_2>0)=\notag\\
 \int_0^\infty\int_0^\infty\int_{-\infty}^{r_1m_1}\int_{-\infty}^{r_2m_2}
f_K(n_1,~ n_2)f_G(m_1,~m_2)dn_2dn_1dm_2dm_1
\end{align}
where $f_K(n_1, n_2) = f_K(D_uw_y({\bf s}), D_uw_y({\bf s}'))$ is the bivariate normal density for the $D_uw_y({\bf s})$ Gaussian
process with parent covariance function $K(\cdot)$ evaluated at ${\bf s}$ and ${\bf s}'$, and similar for $f_G(m_1, m_2)$. Recall
that a univariate Cauchy distribution can be defined as a scale mixture of normals \cite[]{AndrewsMallows1974}. If we write the
corresponding CDF as $F(r)=\int_0^\infty\int_{-\infty}^{rm}\phi(n)\phi(m)dndm$, then the form in Equation \ref{eq:cauch} is
evidently similar and can be regarded as a bivariate analogue. In this way, the spatial Cauchy process is defined such that at
any location the distribution is a univariate Cauchy distribution and for any set of locations the distribution is a sum of
integrals of a similar form.

In the data analysis examples we consider a fixed direction ${\bf u}$ and draw samples of the directional derivative ratio at
each location across a region.  Although we cannot show mean square continuity for this surface, by imposing additional
smoothness conditions on the covariance functions for $w_{y}(\mathbf{s})$ and $w_{x}(\mathbf{s})$, we can argue that this surface
will be almost surely continuous using results from \cite{Kent1989}, following the development in \cite{BanerjeeGelfand2003}.  We
omit the details.

\subsection{Spatial Angular Discrepancy Process}
\label{subsec:disc} At any given location there may be interest not only in the magnitude of the gradients in various directions,
but also in the direction at which the maximum gradients are achieved. Stronger alignment between the directions of maximum
gradients would suggest a stronger relationship between the response and covariate surfaces.

 Consider the covariate process $X({\bf s})$. At location ${\bf s}$ the maximum gradient will be achieved in the direction
described by the unit vector ${\bf u}_X^*=\nabla_X({\bf s})/||\nabla_X({\bf s})||$, and the directional derivative in the
direction of maximum gradient will be $D_{{\bf u}_X^*}X({\bf s}) = ||\nabla_X({\bf s})||$ \cite[]{Banerjee2003}. We can similarly
consider the direction of maximum gradient for the response surface which will occur for ${\bf u}_Y^*=\nabla_Y({\bf
s})/||\nabla_Y({\bf s})||$. In applications, researchers may be interested in the behavior of the response surface, or its
relative rate of change, in the direction ${\bf u}_X^*$.  That is, there may be interest in $D_{{\bf u}_X^*}Y({\bf
s})=\nabla_X({\bf s})'\nabla_Y({\bf s})/||\nabla_X({\bf s})||$ as well as $D_{{\bf u}_X^*}Y({\bf s})/D_{{\bf
      u}_X^*}X({\bf s})=(\nabla_X({\bf s})'\nabla_Y({\bf s}))/(\nabla_X({\bf s})'\nabla_X({\bf s}))$.

At a location ${\bf s}$ the magnitude of the maximum gradient, $||\nabla_X({\bf s})||$, will be the square root of a sum of
squared independent normal random variables; as such, this quantity will have a Chi distribution with $d=2$ degrees of freedom,
possibly scaled by some factor. If the process is isotropic, then $||\nabla_X({\bf s})||$ will have a Chi distribution with $d=2$
degrees of freedom scaled by $SD(D_{{\bf e}_1}X({\bf s}))=SD(D_{{\bf e}_2}X({\bf s}))$. For the Mat\'{e}rn covariance ($\nu=3/2$)
structure this scaling factor will be equal to $\sigma_x\phi_x$.

The unit vector describing the direction of max gradient for the covariate surface can equivalently be described by an angle
$\theta_X({\bf s})$ such that $\tan(\theta_X({\bf s}))=D_{(0,1)}X({\bf s})/D_{(1,0)}X({\bf s})$, and similarly for the response
surface. The angle $\theta$ can take values from $-\pi$ to $\pi$, so inversion of tan must be done with care. Since tan has a
period of only $\pi$, the inverse function is typically taken to be arctan$^*$. \footnote{arctan$^*(S/C)$ is arctan$(S/C)$ if
$C>0$, $S\ge 0$; $\pi/2$ if $C=0$, $S>0$; arctan$(S/C)+\pi$ if $C<0$; arctan$(S/C)+2\pi$ if $C\ge 0$, $S<0$; and undefined if
$C=0$, $S=0$ \cite[]{JammalamadakaSenGupta2001}.}

As smooth functions of a well defined spatial process, the directions of maximum gradient $(\theta_X({\bf s}), \theta_Y({\bf
s}))$ define a bivariate projected Gaussian process, analogous to the univariate projected Gaussian processes described in \cite{WangGelfandSub}. Marginally $\theta_X({\bf s})$ and $\theta_Y({\bf s})$ will each be a projected Gaussian
process \cite[]{WangGelfandSub}. Assuming the Mat\'{e}rn ($\nu=3/2$) covariance structure we can derive the joint distribution of
the two angles at a given location ${\bf s}$:
\begin{align*}
f(\theta_X(\mathbf{s}),\theta_Y(\mathbf{s})) &= \begin{cases}
C\left(\frac{\tilde{A}^{2}(\mathbf{s})\phi(0)}
{ac(ac-\tilde{A}^{2}(\mathbf{s}))} + \frac{\sqrt{2\pi}\tilde{A}(\mathbf{s})}{(ac-\tilde{A}^{2}(\mathbf{s}))^{3/2}}
 L(0,0,\sqrt{\frac{\tilde{A}^{2}(\mathbf{s})}{ac}}) +
 \frac{\phi(0)}{ac}\right),&\hspace{-.05cm} \tilde{A}({\bf s})>0\\
C\left(\frac{\tilde{A}^{2}(\mathbf{s})\phi(0)}
{ac(ac-\tilde{A}^{2}(\mathbf{s}))} + \frac{\sqrt{2\pi}\tilde{A}(\mathbf{s})}{(ac-\tilde{A}^{2}(\mathbf{s}))^{3/2}}
 \left(0.5 - L(0,0,\sqrt{\frac{\tilde{A}^{2}(\mathbf{s})}{ac}})\right) +
 \frac{\phi(0)}{ac}\right),&\hspace{-.05cm}  \tilde{A}({\bf s})<0
\end{cases}
\end{align*}
where $C=\frac{1}{a(2\pi)^{3/2}\sqrt{|\Sigma|}}$, $|\Sigma|=(\sigma_x^2\phi_x^2)^2(\sigma_y^2\phi_y^2)^2$, $a=1/(\sigma_y^2\phi_y^2)$,
$c=(\sigma_y^2\phi_y^2+\beta^2\phi_x^2\sigma_x^2)/(\sigma_x^2\phi_x^2)$, $\tilde{A}(\mathbf{s})=\sqrt{a}\beta
\cos(\theta_X(\mathbf{s})-\theta_Y(\mathbf{s}))$, and $L(0, 0, \rho)$ is the zero mean bivariate normal cdf with correlation
$\rho$ and standard deviations equal to 1 evaluated at $(0,0)$. A
brief derivation is provided in Appendix B. In addition, it is
straightforward to show that at each location ${\bf s}$ the angles, $\theta_X({\bf s})$ and $\theta_Y({\bf s})$, will marginally
be uniform on $(-\pi,\pi)$.

\begin{figure}[h]
\centering
\includegraphics[width=6in]{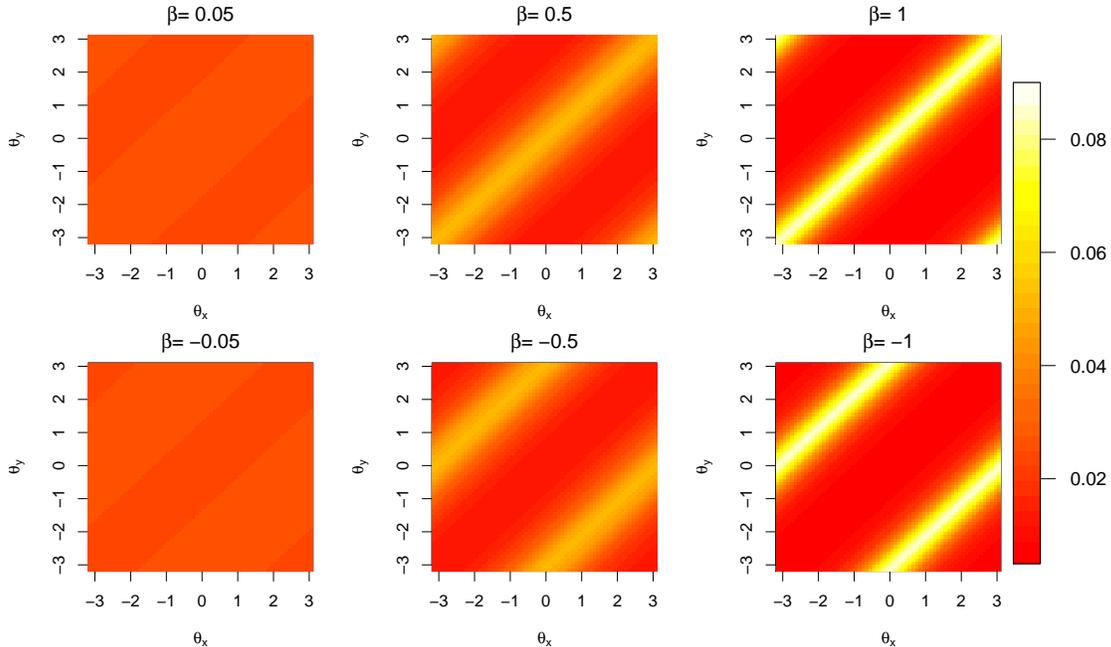}
\caption{Bivariate density for $(\theta_X({\bf s}),\theta_Y({\bf s}))$
  at a fixed ${\bf s}$ for varying values
  of $\beta$. All other parameters are set to the values used for
  simulation in Section \ref{sec:sim}. }
\label{fig:BivThet}
\end{figure}

The above bivariate density is plotted in Figure \ref{fig:BivThet} for $\beta=(\pm 0.05, \pm 0.5, \pm 1)'$, holding all other parameters constant:
$\sigma_y=1=\sigma_x$ and $\phi_y=1.05=\phi_x$. When $\beta>0$ the mass is concentrated around $(\theta_X, \theta_Y)$ pairs that
are equal; when $\beta<0$ the mass is concentrated around pairs where
$\theta_X=\theta_Y-\pi$. When $\beta=\pm 0.05$ the
relationship between $X({\bf s})$ and $Y({\bf s})$ is weak, and the density is roughly uniform over all angle pairs. As the magnitude of $\beta$
increases, the mass becomes increasingly
concentrated around these respective values.

To compare the directions of maximum gradient calculated from the posterior distribution, we compute a ``discrepancy'' for
$\theta_X({\bf s})$ and
$\theta_Y({\bf s})$.
Define: $disc({\bf s}) = 1-\cos(\theta_X({\bf s})-\theta_Y({\bf
  s}))$. As a smooth function of a spatial process, this discrepancy
is a well defined spatial process which we refer to as the {\it spatial angular discrepancy process}. When the maximum gradients
occur in identical directions this process will have a value of zero, when they occur in opposite directions the process will
have a value of two. In our analyses we consider the process across a region, plotting the posterior median surface.

In some areas the $X(\mathbf{s})$ or $Y(\mathbf{s})$ surface may be quite flat and there will be no direction with a gradient
magnitude substantially larger than in the other directions. In these areas a large angular discrepancy may not be as meaningful
as it would be in areas with larger gradient magnitudes. For this reason it is useful to examine these plots in tandem with plots
of the local directional sensitivity process described in Section \ref{subsec:direcderiv}. 

\subsection{Extensions to Non-Gaussian Data Models Using the Chain Rule}
\label{subsec:chain} One can imagine modeling scenarios where the response will be non-Gaussian, for instance, binary or
zero-inflated data. In these cases, a latent Gaussian surface may be utilized and spatial gradients can still be considered. The
resulting gradient surfaces will be for the latent Gaussian process, not necessarily the response of interest; however,
differentiable transformations will allow for inference on the non-Gaussian response surface. Recall the chain rule result
presented in \cite{Majumdar2006}: for $g(\cdot)$ differentiable on $\mathbb{R}^1$ and $W({\bf s})=g(V({\bf s}))$, the directional
derivative $D_{\bf u}W({\bf s})$ exists and is given by $D_{\bf u}W({\bf s}) = D_{\bf u}g(V({\bf
  s}))=g'(V({\bf s}))D_{\bf u}V({\bf s})$. Here $V({\bf s})$ is the latent Gaussian surface and the function
$g(\cdot)$ describes its relationship to the non-Gaussian response of interest.%

When analyzing the
Duke Forest data in Section \ref{sec:data}, we work with a log-Gaussian Cox process model for the intensity surface. To study
gradients associated with the intensity surface itself, we can use the fact that, if $Y({\bf s})$ is our response and $Z({\bf
s})=\log Y({\bf s})$, then $D_{\bf u}Y({\bf s}) = \exp(Z({\bf s}))D_{\bf u}Z({\bf s})$.

In the case of binary data we can write a model with response $Y({\bf
  s})=1$ if $Z({\bf s})>0$ and $Y({\bf s})=0$ if $Z({\bf s})<0$,
where $Z({\bf s})$ is a latent Gaussian process centered at a linear function of a spatial covariate $X({\bf s})$. We have a
spatial probit regression since $P(Y({\bf s})=1)=\Phi(Z({\bf s}))$ \cite[]{HeagertyLele1998}. There may be interest in identifying
areas where the transition between regions of high absence probabilities to regions of low presence probabilities is rapid (or
vice versa). Consideration of the local directional sensitivity process associated with the probability surface $P({\bf s})
\equiv P(Y({\bf s})=1)$ and the covariate surface $X({\bf s})$ would provide an avenue to answer these kinds of questions. Now,
with the function $g(\cdot)=\Phi(\cdot)$, we obtain $D_{\mathbf{u}}P(\mathbf{s})=
\phi(Z(\mathbf{s}))D_{\mathbf{u}}Z(\mathbf{s})$.

\section{Simulation Example}
\label{sec:sim} We consider a simulation example to explore the behavior of gradient quantities in a controlled setting. From the
foregoing development, in the context of spatial gradient analysis the quantities of interest will be directional derivatives.
These derivatives are unobservable even in a simulation study with known parameters. Thus, assessing inference performance with
regard to these quantities requires some novelty.

Recall that the gradient processes describe the shape and behavior of spatial surfaces. With a simulated dataset, we are able to
draw a realization of the Gaussian process over a larger number of locations, allowing for fairly detailed understanding of the
spatial surfaces. To avoid any unfair advantage that might come from the increased sample size, we use only a subset of the
locations to fit the model and reserve the full set of locations for assessing the quality of our inference. Contour lines
highlighting the shape of the surface are interpolated using the full set of locations. Conclusions made using gradients are then
compared to those suggested by the contour lines to examine performance.

\subsection{Data and Model}
We simulate $X(\mathbf{s})$ a realization from a mean zero Gaussian process on $[0,10]\times[0,10]$, and
$Y(\mathbf{s})|X(\mathbf{s})$ a realization from a Gaussian process with mean $\beta X(\mathbf{s})$. We use Mat\'{e}rn
covariance functions setting $\nu=3/2$ in order to capitalize on the resultant closed form. We simulate the Gaussian processes
assuming Mat\'{e}rn covariance functions with parameters $\phi_x=1.05=\phi_y$, $\sigma^2_x=1=\sigma^2_y$, $\beta_0=0=\alpha_0$,
and $\beta_1=0.5$. We draw a larger realization at 2000 locations, to allow for finer knowledge of the underlying surface, from
which we consider a subset of 200 locations to be our ``observations''.

Treating $(Y({\bf s}),X({\bf s}))'$ as a multivariate Gaussian process, we fit a coregionalization model using the conditional
parameterization, as described in Section \ref{sec:dists}. The model is fitted using the \texttt{`spbayes'} package in \texttt{R}
by first fitting parameters for $X({\bf s})$, then fitting parameters for $Y({\bf s})|X({\bf s})$.  We obtain 2000 samples after a burn-in of 500 iterations
and a thinning of every fifth iterate. We assume $\nu=3/2$. Priors for the remaining parameters are: $\alpha_0  \sim N(0,100),
\beta_0  \sim N(0,100), \beta_1 \sim N(0,100), \phi_x, \phi_y \sim U(0.5, 10), \sigma^2_x, \sigma^2_y \sim IG(2,0.1)$
where $\alpha_0=E(X({\bf s}))$ and $\beta_0+\beta_1 X({\bf s}) = E(Y({\bf s})|X({\bf s}))$. Summaries of the posterior parameter samples are provided in Table \ref{tab:XYGP}.

\begin{figure}[h]
\centering
\includegraphics[width=4.5in]{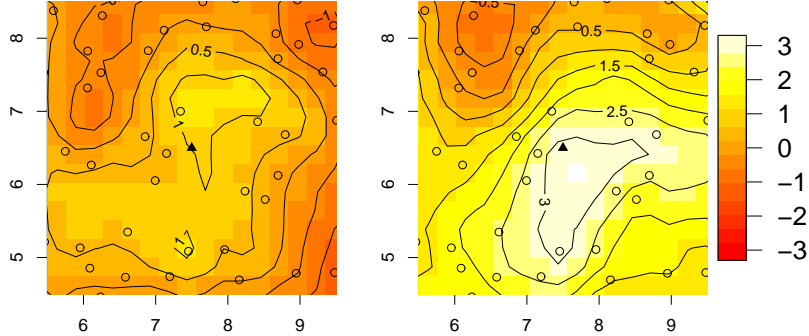}
\caption{$X({\bf s})$ (left) and $Y({\bf s})$ (right) subregions around
  ${\bf s}^*$,
  where we estimate the gradient.} 
\label{fig:loc.st}
\end{figure}

\begin{table}[ht]
\begin{center}
\begin{tabular}{r|ccc|c}
  \hline
Parameter &0.025 & Mean & 0.975 & Truth\\
  \hline
$\alpha_0$ & -0.7126 & -0.0614 & 0.6214 & 0\\
$\beta_0$ & 0.2726 & 0.8179 & 1.3867 & 0 \\
$\beta_1$ & 0.4685 & 0.5943 & 0.7202 & 0.5\\
$\sigma^2_x$ & 0.6725 & 1.0692 & 1.7390 & 1\\
$\phi_x$ & 0.8242 & 1.0230 & 1.2124 & 1.05\\
$\sigma^2_y$ & 0.4911 & 0.8057 & 1.3341 & 1\\
$\phi_y$ & 0.8572 & 1.0718 & 1.3116 & 1.05\\
   \hline
\end{tabular}
\end{center}
\caption{Parameter estimates for the $X({\bf s})$ and $Y({\bf s})|X({\bf s})$ model.}
\label{tab:XYGP}
\end{table}
We consider a region centered at the location ${\bf s}^*=(7.5,6.5)$. The $X({\bf s})$ and $Y({\bf s})$ values at this location
are provided in Figure \ref{fig:loc.st}. (The full processes were realized on $[0,10]\times[0,10]$, but we only show a subregion
here.) The interpolated surface and contour lines are produced using the full 2000 locations, while the circles indicate the
subset of
200 locations used to predict the gradient. 

\subsection{Local Directional Sensitivity Process}
We are interested in the behavior of $D_{\bf u}Y({\bf s})/D_{\bf
  u}X({\bf s})$. We consider ${\bf u}=(1,0)$ and ${\bf
  u}=(0,1)$. Since $D_{-{\bf u}}Y({\bf s})=-D_{\bf u}Y({\bf s})$
\cite[]{Banerjee2003}, any discussion of the behavior in ${\bf u}$ direction implies
the opposite
behavior is occurring in the opposite direction. When applied to the ratios, this means that the local
directional sensitivity process will be equal for ${\bf u}$ and $-{\bf u}$.

Returning to the region in Figure \ref{fig:loc.st}, to visualize the local directional sensitivity process we draw samples at a
grid of 125 locations denoted as $\{{\bf
  s}^*_1,\hdots,{\bf s}^*_{n}\}$. We draw 2000 samples of
$(\nabla_Y({\bf s}^*_1),\hdots,\nabla_Y({\bf s}^*_n), \nabla_X({\bf
  s}^*_1),\hdots, \nabla_X({\bf s}^*_n))'$ from the joint predictive distribution, again, given
the fitting data at 200 observed locations. For each of these samples
we calculate $D_{\bf u}Y({\bf s}^*)/D_{\bf u}X({\bf s}^*)$. We
summarize the central behavior of these Cauchy quantities at a given location using the median value of the ratios.

\begin{figure}[h]
\centering
\includegraphics[width=5in]{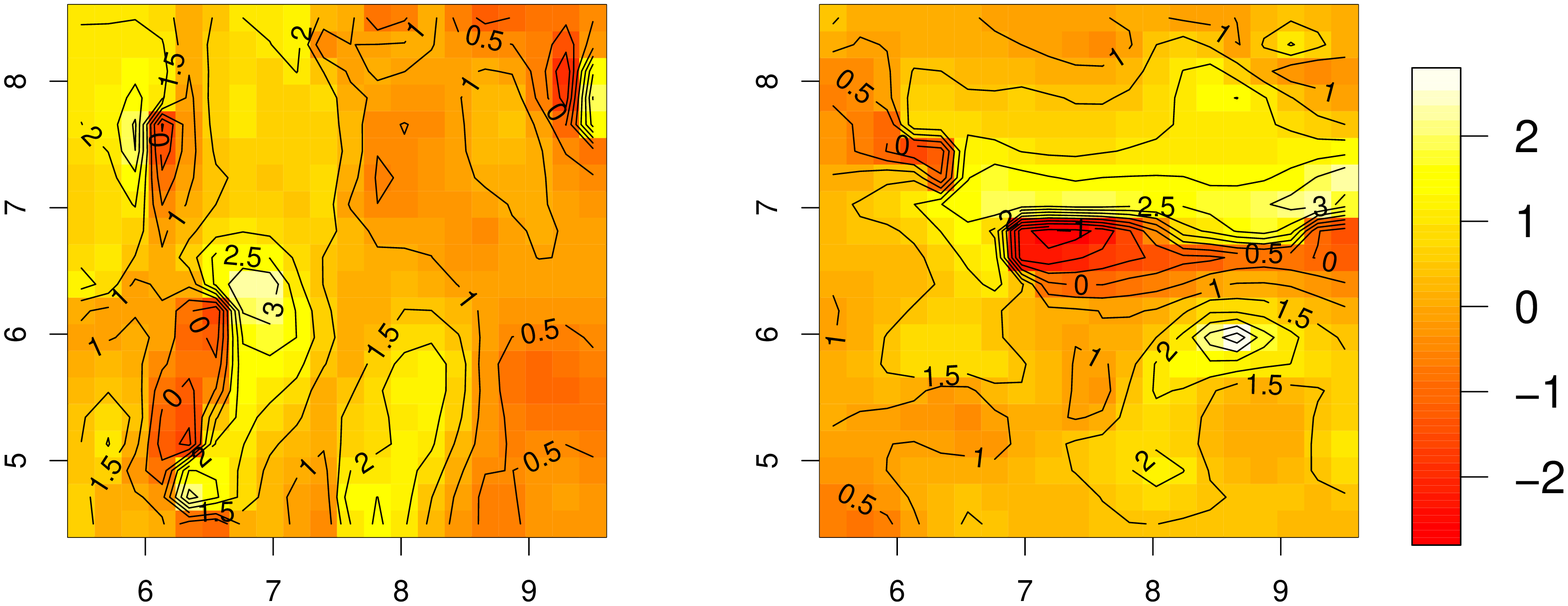}
\caption{Posterior median of $D_{\bf u}Y({\bf s})/D_{\bf u}X({\bf s})$ in the directions ${\bf u}=(1,0)$ (left) and ${\bf u}=(0,1)$ (right). }
\label{fig:RatDeriv_posterior}
\end{figure}

For the two directions being considered, we plot the median predictive surface in Figure \ref{fig:RatDeriv_posterior}.
Interpretation of these surfaces requires examination of the sign of the ratio as well as the magnitude. Magnitudes less than 1
suggest that the $X({\bf s})$ surface is changing more rapidly than the $Y({\bf s})$ surface; magnitudes greater than 1 suggest
that the $Y({\bf s})$ surface is
 changing more rapidly than the $X({\bf s})$ surface; negative values suggest
 that one surface is increasing while the other decreases; positive
 values suggest that both surfaces are either increasing or decreasing.

The direction ${\bf u}=(1,0)$ points towards the east. The corresponding ratio surface is provided in the left hand plot of
Figure \ref{fig:RatDeriv_posterior}. There is a peak in the ratio surface around (6.75,6.25), suggesting that both surfaces are
either decreasing or increasing and that the $Y({\bf s})$ surface is doing so more rapidly. Referring back to Figure
\ref{fig:loc.st}, the contour lines indicate that both surfaces are increasing at that location looking east, and that the
$Y({\bf s})$ surface is doing so more rapidly. The direction ${\bf u}=(0,1)$ points towards the north. The corresponding ratio
surface is provided in the right hand plot of Figure \ref{fig:RatDeriv_posterior}. There is a peak in the ratio surface around
(8.75, 6), suggesting that both surfaces are either decreasing or increasing and that the $Y({\bf s})$ surface is doing so more
rapidly. Referring back to Figure \ref{fig:loc.st}, the contour lines indicate that both surfaces are increasing at that location
looking east, and that the $Y({\bf s})$ surface is doing so more rapidly.

\subsection{Spatial Angular Discrepancy Process}
Consider again the region in Figure \ref{fig:loc.st} and the grid of 125 locations denoted as $\{{\bf
  s}^*_1,\hdots,{\bf s}^*_{n}\}$. Sampling gradients from the joint predictive distribution for
$(\nabla_X({\bf s}_i^*), \nabla_Y({\bf s}_i^*))'$, we calculate the direction of maximum gradient as $\nabla_X({\bf
s}_i^*)/||\nabla_X({\bf s}_i^*)||$ and $\nabla_Y({\bf s}_i^*)/||\nabla_Y({\bf s}_i^*)||$ for each sample gradient at each
location in the figure. Denote these angles (in radians) as
$\theta_X({\bf s}_i^*)$ and $\theta_Y({\bf s}_i^*)$ respectively. 

\begin{figure}[h]
\centering
\includegraphics[width=2.5in]{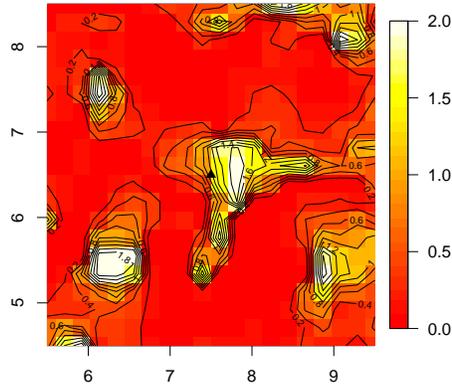}
\caption{Posterior median $disc({\bf s})$.}
\label{fig:distMax}
\end{figure}

 We compute the discrepancy
between these angles at each location as in Section \ref{subsec:disc} and provide the posterior median values in Figure \ref{fig:distMax}. Most
of the region has an associated distance of 0, suggesting that both $X(\mathbf{s})$ and $Y(\mathbf{s})$ are typically increasing
most rapidly in the same direction. However, there are a few small regions where the distance peaks towards a value of 2,
locations where the $X(\mathbf{s})$ and $Y(\mathbf{s})$ surfaces are increasing
in nearly opposite directions.

\section{Duke Forest Point Pattern, Elevation Example}
\label{sec:data}
\subsection{Data}
Our illustrative data set is a collection of point patterns of tree species present at the Blackwood site in the Duke Forest in
Durham, NC. The site is 5 hectares in area and exhibits a range of elevation. A road and powerline separate the site into three
subregions, and we focus on the northwestern of these subregions. We consider two tree species, Flowering Dogwood ({\it Cornus
florida}) and Sweetgum ({\it Liquidambar styraciflua}), with regard to their respective point patterns of locations within the
site in the year 2000. Some trees have multiple stems observed at a single location; however we treat these as a single observed
tree at the given location.  The point pattern for Sweetgum consists of 531 trees, and the point pattern for Flowering Dogwood
consists of 570 trees. Elevation is recorded at each location where a tree of any species was observed, resulting in 5654
elevation observations.

\begin{figure}[h]
\centering
\includegraphics[width=6.5in]{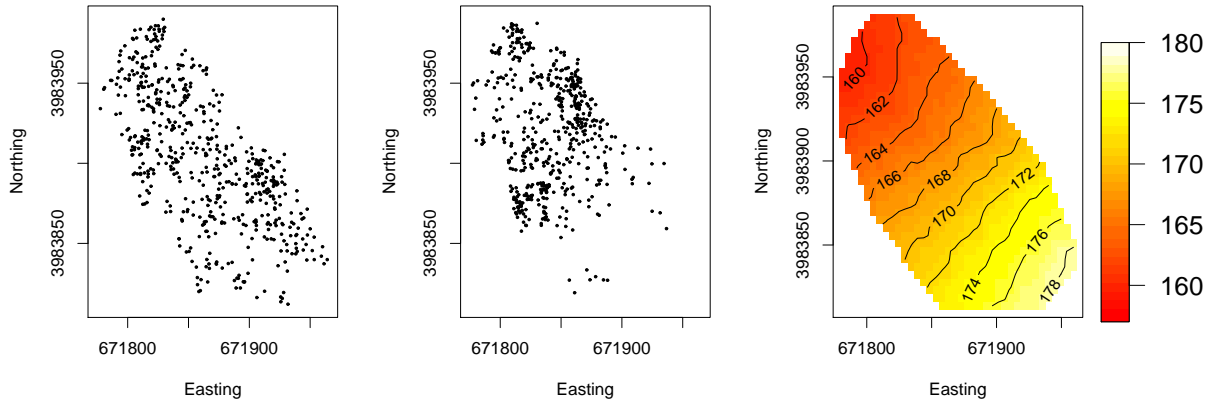}
\caption{From left to right, observed point patterns for Flowering Dogwood and Sweetgum and observed elevation.}
\label{fig:elevpointpattern}
\end{figure}

Figure \ref{fig:elevpointpattern} provides a heatmap of the elevation data. There is a clear increase in elevation across the
region in a roughly southeastern direction. Figure \ref{fig:elevpointpattern} also provides the observed point patterns for each
of the species. Flowering Dogwood is well dispersed across the entire region while Sweetgum is more abundant in
the northwestern half of the region.

\subsection{Point Pattern Model given Elevation}
Given the observed elevation, $X({\bf s})$, we model the intensity for each species as a log-Gaussian Cox process: $\lambda({\bf
s})  = \exp(\beta_0 + \beta_1X({\bf s})) \lambda_0({\bf s})$, $\lambda_0({\bf s})  = \exp(w_z({\bf s}))$, and  $w_z({\bf s}) \sim
GP(0, \rho(\cdot|\phi_z, \sigma_z^2))$. Placing priors on the parameters, we have the following model:
\begin{align*}
[\tilde{Y}|\beta_0, \beta_1, \{\lambda_0({\bf s}), {\bf s}\in D\}][\lambda_{0}({\bf s})|\sigma^2, \phi][\beta_0, \beta_1][\phi][\sigma^2]\\
 = [\tilde{Y}|\beta_0, \beta_1, \{w_z({\bf s}), {\bf s}\in D\}][w_z({\bf s})|\sigma^2, \phi][\beta_0, \beta_1][\phi][\sigma^2]
\end{align*}
where $\tilde{Y}$ is the observed point pattern and $\lambda_0({\bf s})$ is equivalently considered using $w_z({\bf s})$.

We approximate the likelihood by dividing the region into a fine grid with cells $\{A_l; l=1,\hdots, L\}$. This gives us the likelihood
\begin{align*}
L(\lambda({\bf s}), {\bf s}\in D; {\bf s}_1, {\bf s}_2, \hdots, {\bf s}_n) & \approx \Pi_i \lambda({\bf s}_i)\exp(-\lambda(D))\\
\lambda(D) & \approx \Sigma_l \exp(X'(A_l)\beta + w_z(A_l))
\end{align*}
where $w_z(A_l)$ corresponds to a realization from a Gaussian process evaluated at a
representative point in each grid cell $A_l$. This likelihood can be
sampled using elliptical slice sampling, as described in \cite{Murray2009} and \cite{MurrayAdams2010}.

\subsection{Conditional Bivariate Model}
In terms of $Z(\mathbf{s}) = \log(\lambda(\mathbf{s}))$, we immediately have a conditional bivariate GP model with elevation:
$Z({\bf s})|X({\bf s})  = \beta_0 + \beta_1X({\bf s}) + w_z({\bf s})$ and $X({\bf s}) = \alpha_0 + w_x({\bf s})$ with $w_z({\bf
s}) \sim GP(0,\rho(\cdot|\sigma^2_z, \phi_z))$ and $ w_x({\bf s}) \sim GP(0,\rho(\cdot|\sigma^2_x, \phi_x))$. Hence, we are in
the framework developed above and can apply the proposed gradient analyses.  The difference in this case will be that the
Gaussian response $Z({\bf s})$ is latent and thus unobserved. The uncertainty about $Z({\bf s})$ is propagated through the model
by drawing a posterior sample of the $Z({\bf s})$ surface for each posterior sample of the parameters.

\begin{figure}[h]
\centering
\includegraphics[width=5in]{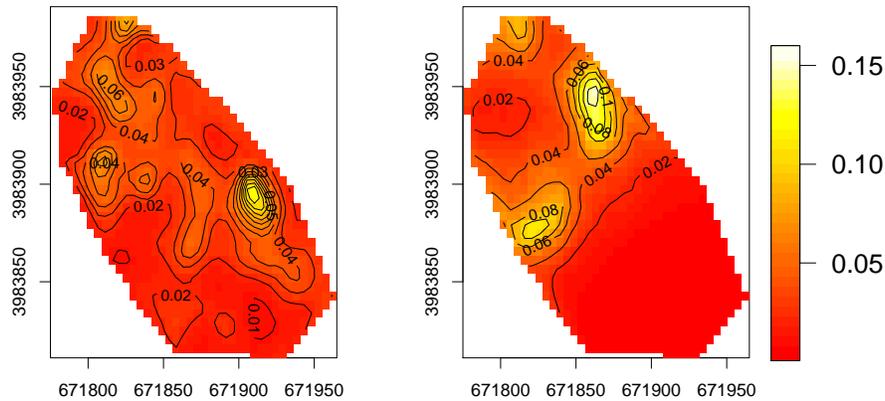}
\caption{Posterior median of the intensity surface for Flowering
  Dogwood (left) and Sweetgum (right).}
\label{fig:Lam}
\end{figure}

Table \ref{tab:GPelev} provides the fitted parameter values for elevation. Tables \ref{tab:GPcofl} and \ref{tab:GPlist} provide
the fitted parameter values for each of the species. 
Note that the $\phi_z$ parameter is fixed at the minimum contrast
estimate, as suggested in \cite{Moller1998}, to facilitate identifiability in the fitting of the
log-Gaussian Cox process. The fixed
values are $\phi_z=0.1063$ and 0.0434 for Flowering Dogwood and
Sweetgum respectively.
The 95\% credible interval for $\beta_1$ contains zero for Flowering Dogwood, but not for Sweetgum.
This negative coefficient suggests that the intensity of Sweetgum
decreases as elevation increases, while the intensity of Flowering
Dogwood is not responsive to elevation changes at this scale.  Figure \ref{fig:Lam} provides the posterior median intensity for
each of the species. Both intensities have fairly low values across most of the domain, with a few regions of higher intensity.

\begin{table}[ht]
\begin{center}
\begin{tabular}{r|ccc}
  \hline
Parameter &0.025 & Mean & 0.975\\
  \hline
$\alpha_0$ & 166.1043 & 167.8577 & 169.4328\\
$\sigma^2_x$ & 7.1681 & 9.4117 & 12.9494\\
$\phi_x$ & 0.0791 & 0.0892 & 0.0979\\
   \hline
\end{tabular}
\end{center}
\caption{Parameter estimates for elevation.}
\label{tab:GPelev}
\end{table}

\subsection{Local Directional Sensitivity Process}

\begin{table}[ht]
\begin{center}
\begin{tabular}{r|ccc}
  \hline
Parameter &0.025 & Mean & 0.975\\
  \hline
$\beta_0$ & -3.6416 & -3.5478 & -3.4590\\
$\beta_1$ & -0.0866 & -0.0284 & 0.0348\\
$\sigma^2_z$ & 0.3367 & 0.5237 & 0.7648\\
   \hline
\end{tabular}
\end{center}
\caption{Parameter estimates for $Z({\bf s})|X({\bf s})$ model for Flowering Dogwood.}
\label{tab:GPcofl}
\end{table}

\begin{table}[ht]
\begin{center}
\begin{tabular}{r|ccc}
  \hline
Parameter &0.025 & Mean & 0.975\\
  \hline
$\beta_0$ & -4.2161 & -4.0297 & -3.8728\\
$\beta_1$ & -0.4094 & -0.2619 & -0.1024\\
$\sigma^2_z$ & 0.7942 & 1.2652 & 1.9160\\
   \hline
\end{tabular}
\end{center}
\caption{Parameter estimates for $Z({\bf s})|X({\bf s})$ model for Sweetgum.}
\label{tab:GPlist}
\end{table}

The methods developed in previous sections allow for straight forward examination of $D_{\bf u}\texttt{log}\lambda({\bf
s})/D_{\bf u}X({\bf s})= D_{\bf u}Z({\bf s})/D_{\bf u}X({\bf s})$, although interest is more likely in the behavior of the
intensity surface itself. Applying the spatial gradient chain rule discussed in Section \ref{subsec:chain} gives us $D_{\bf
u}\lambda({\bf s})=\exp(Z({\bf s}))D_{\bf u}Z({\bf s})$, with the directional derivative ratio $D_{\bf u}\lambda({\bf s})/D_{\bf
u}X({\bf s})=\exp(Z({\bf s}))D_{\bf u}Z({\bf s})/D_{\bf u}X({\bf s})$. As before, we can simplify this in terms of the
independent Gaussian processes:
\begin{align*}
\frac{D_{\bf u}\lambda({\bf s})}{D_{\bf u}X({\bf s})} & = \exp(Z({\bf s}))\frac{D_{\bf u}Z({\bf s})}{D_{\bf u}X({\bf s})}\\
& = \frac{\exp(Z({\bf s}))[\beta_1D_{\bf u}w_x({\bf s}) + D_{\bf u}w_z({\bf s})]}{D_{\bf u}w_x({\bf s})} \\
& = \exp(Z({\bf s}))[\beta_1 + \frac{D_{\bf u}w_z({\bf s})}{D_{\bf u}w_x({\bf s})}]
\end{align*}
As before, $\beta_1$ will center the Cauchy random variable, but now there will also be scaling according to the value of
$\exp(Z({\bf s}))$.

\begin{figure}[ht]
\centering
\includegraphics[width=5in]{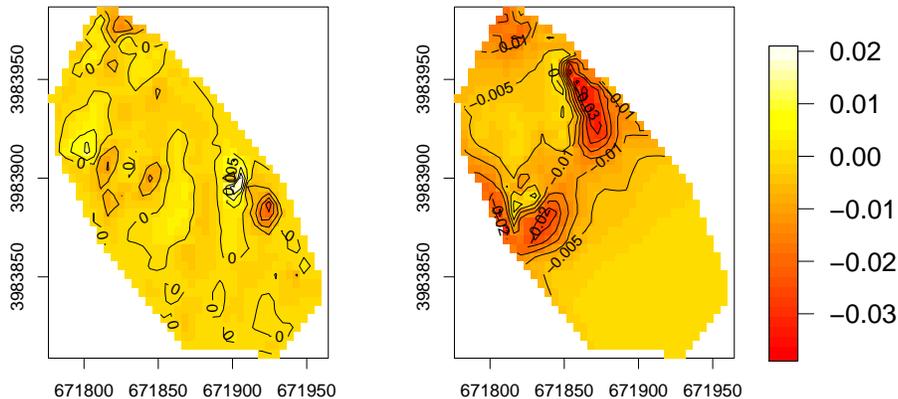}
\caption{Posterior median of $D_{\bf u}\lambda({\bf s})/D_{\bf u}X({\bf s})$ for
  Flowering Dogwood (left) and Sweetgum (right);
${\bf u}=(0.8508, -0.5255)$.}
\label{fig:cofllistDYDXuu}
\end{figure}

In Figure \ref{fig:elevpointpattern} we saw a clear increase in elevation in a roughly southeastern direction. We approximate
this direction by the unit vector ${\bf u}=(0.8508, -0.5255)$. Again, via the chain rule, we can consider the behavior of the
directional derivative ratios in this direction for each of the species. Figure \ref{fig:cofllistDYDXuu} plots the resulting
posterior median $D_{\bf u}\lambda({\bf s})/D_{\bf u}X({\bf s})$ surfaces.



For Flowering Dogwood the majority of the domain has a ratio close to zero. This suggests that the changes in the intensity are
negligible compared to the changes in elevation. Recalling the fairly even spread of the trees in the region, as well as the
non-significant $\beta_1$, this pattern makes sense.

For Sweetgum, virtually the entire region has a negative directional derivative ratio. This aligns with our interpretation of the
significantly negative $\beta_1$ coefficient, namely that as elevation increases the intensity decreases. There are a few
subregions where the change in intensity occurs more rapidly than elsewhere, and there is a larger subregion where the change in
intensity is zero due to an absence of trees. The cause for subregions of rapid change could be further illuminated through
examination of other factors in those regions. Similarly, the region of zero change in intensity could be roughly interpreted as
Sweetgum having an aversion to elevations beyond a certain value.


\subsection{Spatial Angular Discrepancy Process}
\begin{figure}[tbp]
\centering
\includegraphics[width=5in]{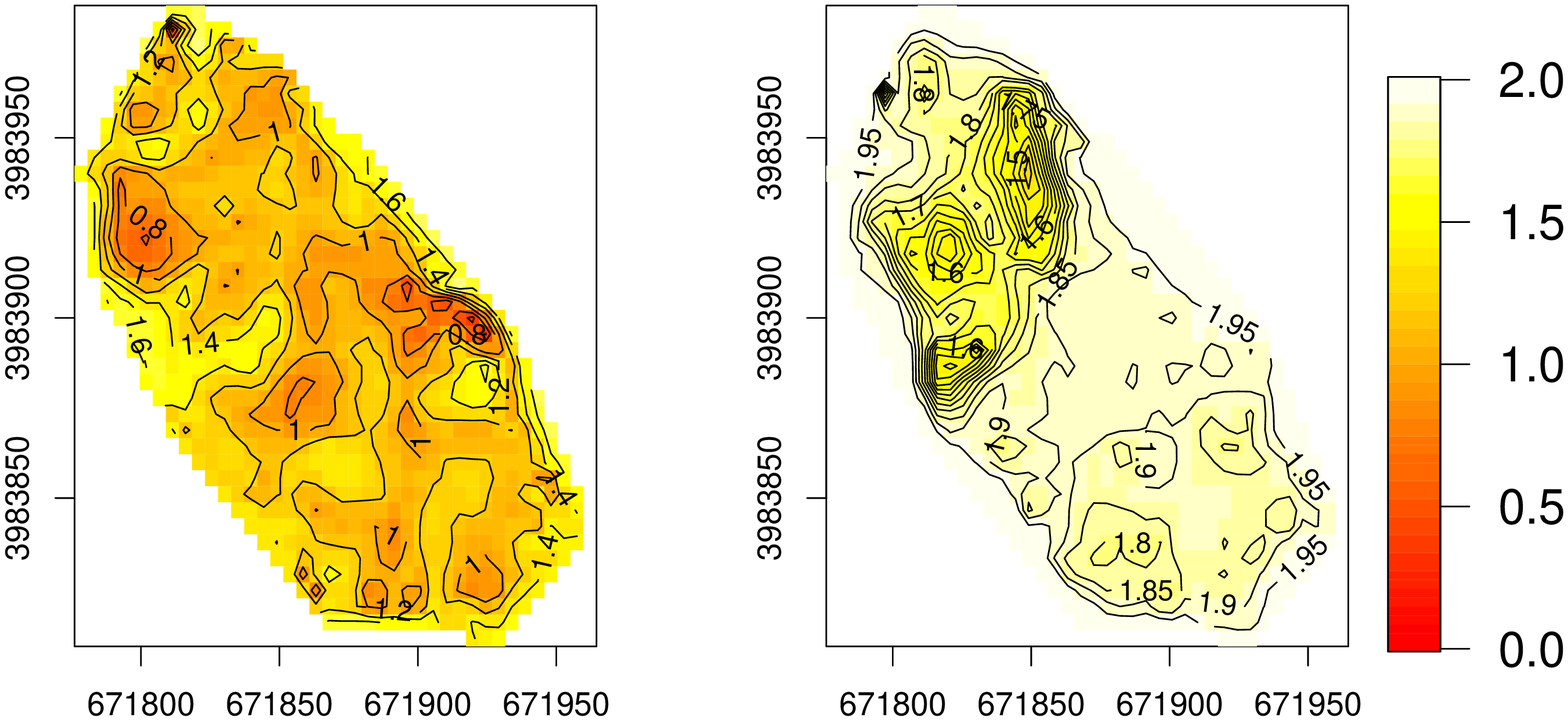}
\caption{Posterior median $disc({\bf s})$ for Flowering Dogwood
  (left) and
  Sweetgum (right).}
\label{fig:distcofllist}
\end{figure}

Finally, we can compare the intensity and elevation surfaces by computing the discrepancy between their directions of maximum
gradient at each location, i.e., the posterior median of $disc({\bf s})=1-\cos(\theta_X({\bf s})-\theta_Y({\bf s}))$ across the
region. Values close to 2 suggest the surfaces are most rapidly increasing in opposite directions; values close to 0 suggest the
surfaces are most rapidly increasing in the same direction. The posterior median discrepancy surfaces are provided in Figure
\ref{fig:distcofllist} for Flowering Dogwood and Sweetgum.

The discrepancies for Flowering Dogwood roughly range between 0.8 and
1.5. There is no
clear pattern, which supports there being no strong relationship
between Flowering Dogwood intensity and elevation.

The pattern for Sweetgum is quite different. All of the discrepancies appear to be between 1.5 and 2, with most around 1.9. This
suggests that the Sweetgum intensity and elevation are increasing in nearly opposite directions virtually everywhere in the
domain. This again confirms the negative relationship, and additionally highlights this pattern as being slightly weaker in the
northern part of the region.

\section{Summary and Conclusions}
\label{sec:summary} We have developed methodology for performing
spatial sensitivity analysis for a bivariate process where one
variable is treated as a response to the other variable. Consideration
of the associated directional derivatives can be done jointly and
results in a multivariate Gaussian process directly derivable from the
model for the parent process fit in a Bayesian framework. Utilizing
the posterior draws of the process parameters, all gradient analysis
occurs post model fitting. 

Using the directional derivatives, we proposed two derived processes
in order to learn about the relationship between the
response and covariate processes. The first is the local directional
sensitivity process, inspired by quantities explored in standard
sensitivity analysis. This process captures local variation in the
relationship between the two variables and provides deeper insight
into their relative behavior across the region. The second is the
spatial angular discrepancy process, capturing the discrepancy between
the directions in which the process surfaces are most rapidly
increasing. Spatial plots of this discrepancy surface highlight
regions of the domain where the two processes behave most similarly
and most differently, again informing on the variables' spatial relationship. 

Our application involved elevation data and point patterns of trees
collected from the Duke Forest. Using a log-Gaussian Cox process, we
studied local directional sensitivity of the intensity to elevation
through the use of a spatial gradient chain rule. Through the two tree
species, this example illustrates the different results one would
expect from a spatial sensitivity analysis when the variables are
significantly related versus when they are not.

The current theory provides opportunity for several extensions and
applications. Many ecological data sets are observed at multiple time
points, in part to see if the relationships between the variables of
interest are changing over  time. With this in mind, future work on
gradient analyses may involve the incorporation of temporal
effects. Ecological data sets can also have multiple responses
such as leaf traits that are being related to multiple climate
covariates \cite[e.g.][]{Thuiller2004}, any or all of which may have
relationships that 
could be better highlighted through a spatial gradient analysis under a joint
model. Finally, the spatial gradient chain rule will allow us to
consider novel non-Gaussian responses modeled through latent Gaussian
process models as described in Section \ref{subsec:chain}. 
\appendix
\section{}
We show that $D_{\bf u}w_y({\bf s})/D_{\bf
    u}w_x({\bf s})$ is a well defined spatial stochastic process.

First consider two locations ${\bf s}$ and ${\bf s}'$. For notational convenience,
let $D_{\bf u}w_y({\bf s})=n_1$, $D_{\bf u}w_y({\bf s}')=n_2$,  $D_{\bf u}w_x({\bf s})=m_1$, and $D_{\bf
  u}w_x({\bf s}')=m_2$. We can write the joint distribution of the ratios in
terms of the joint distributions for the two Gaussian processes:
\begin{align*}
P(\frac{n_1}{m_1}<r_1,~ \frac{n_2}{m_2}<r_2)
& = P(n_1<r_1m_1,~ n_2<r_2m_2,~ m_1>0,~ m_2>0)\\
& +P(n_1<r_1m_1,~ n_2>r_2m_2,~ m_1>0,~ m_2<0)\\
& + P(n_1>r_1m_1~, n_2<r_2m_2,~ m_1<0,~ m_2>0)\\
&+ P(n_1>r_1m_1,~ n_2>r_2m_2,~ m_1<0,~ m_2<0)\\
& = \int_0^\infty\int_0^\infty\int_{-\infty}^{r_1m_1}\int_{-\infty}^{r_2m_2}
f_K(n_1, n_2)f_G(m_1,m_2)dn_2dn_1dm_2dm_1 \\
& +\int_0^\infty\int_{-\infty}^0\int_{-\infty}^{r_1m_1}\int_{r_2m_2}^{\infty}
f_K(n_1, n_2)f_G(m_1,m_2)dn_2dn_1dm_2dm_1 \\
& +\int_{-\infty}^0\int_0^\infty\int_{r_1m_1}^{\infty}\int_{-\infty}^{r_2m_2}
f_K(n_1, n_2)f_G(m_1,m_2)dn_2dn_1dm_2dm_1 \\
& +\int_{-\infty}^0\int_{-\infty}^0\int_{r_1m_1}^{\infty}\int_{r_2m_2}^{\infty}
f_K(n_1, n_2)f_G(m_1,m_2)dn_2dn_1dm_2dm_1
\end{align*}
The above joint distribution depends only on integrals of multivariate
normal densities, so consistency across permutations of the labels is
clearly satisfied.

At a single location the directional derivative ratio is a
ratio of two independent normal random variables. As such, it will
have a Cauchy distribution with a scale parameter dependent on the
parameters in the covariance functions $K(\cdot)$ and $G(\cdot)$,
namely $SD(n_1)/SD(m_1)$. Next
we show that marginalizing the bivariate distribution over $r_2$ will
reduce to the known distribution for the univariate case:
\begin{align*}
P(\frac{n_1}{m_1}<r_1) & =
\int_{-\infty}^\infty P(\frac{n_1}{m_1}<r_1,~ \frac{n_2}{m_2}<r_2) dr_2 \\
& =
\int_{-\infty}^\infty\int_0^\infty\int_0^\infty\int_{-\infty}^{r_1m_1}\int_{-\infty}^{r_2m_2}
f_K(n_1, n_2)f_G(m_1,m_2)dn_2dn_1dm_2dm_1dr_2 \\
& +\int_{-\infty}^\infty\int_{-\infty}^0\int_0^\infty\int_{-\infty}^{r_1m_1}\int_{r_2m_2}^{\infty}
f_K(n_1, n_2)f_G(m_1,m_2)dn_2dn_1dm_2dm_1dr_2 \\
& +\int_{-\infty}^\infty\int_{-\infty}^0\int_0^\infty\int_{r_1m_1}^{\infty}\int_{-\infty}^{r_2m_2}
f_K(n_1, n_2)f_G(m_1,m_2)dn_2dn_1dm_2dm_1dr_2 \\
& +\int_{-\infty}^\infty\int_{-\infty}^0\int_{-\infty}^0\int_{r_1m_1}^{\infty}\int_{r_2m_2}^{\infty}
f_K(n_1, n_2)f_G(m_1,m_2)dn_2dn_1dm_2dm_1dr_2 \\
& =
\int_0^\infty\int_0^\infty\int_{-\infty}^{r_1m_1}f_K(n_1)f_G(m_1,m_2)dn_1dm_2dm_1\\
&+\int_0^\infty\int_{-\infty}^0\int_{-\infty}^{r_1m_1}f_K(n_1)f_G(m_1,m_2)dn_1dm_2dm_1\\
&+\int_{-\infty}^0\int_{0}^\infty\int_{r_1m_1}^{\infty}f_K(n_1)f_G(m_1,m_2)dn_1dm_2dm_1\\
&+\int_{-\infty}^0\int_{-\infty}^0\int_{r_1m_1}^{\infty}f_K(n_1)f_G(m_1,m_2)dn_1dm_2dm_1\\
& = \int_0^\infty\int_{-\infty}^{r_1m_1}f_k(n_1)f_G(m_1)dn_1dm_1 + \int_{-\infty}^0\int_{r_1m_1}^\infty f_k(n_1)f_G(m_1)dn_1dm_1
\end{align*}
We can then rewrite this as:
\begin{align*}
P(\frac{n_1}{m_1}<r_1) & = \int_0^\infty
F_K(r_1m_1)f_G(m_1)dm_1+\int_{-\infty}^0 F_K(-r_1m_1)f_G(m_1)dm_1\\
& = 2\int_0^\infty F_K(r_1m_1)f_G(m_1)dm_1
\end{align*}
The associated density will then be $2\int_0^\infty
m_1f_K(r_1m_1)f_G(m_1)dm_1$, which was shown in
\cite{AndrewsMallows1974} to be a Cauchy distribution. In this case we
will have a Cauchy with scale parameter $SD(n_1)/SD(m_1)$, as
desired.

The marginalization is straight forward for larger dimensions with
marginalization over the two normal densities occurring in a similar way.

\section{}
We are interested in the marginal behavior of $f(\theta_X({\bf s}),
\theta_Y({\bf s}))$ at a location ${\bf s}$ where
$\tan(\theta_X({\bf s}))=D_{(0,1)}X({\bf s})/D_{(1,0)}X({\bf s})$,
$\tan(\theta_Y({\bf s}))=D_{(0,1)}Y({\bf s})/D_{(1,0)}Y({\bf s})$, and\\
$(\nabla_X({\bf s}), \nabla_Y({\bf s}))'= (D_{(1,0)}Y({\bf s}), D_{(0,1)}Y({\bf s}),
D_{(1,0)}X({\bf s}),D_{(0,1)}X({\bf s}))'$ are
normally distributed with the covariance structure provided in
Equation 1.

After converting to polar coordinates, the change of variables formula
gives the following integral, where $g(\cdot)$ is the multivariate normal
density associated with $(\nabla_Y({\bf s}), \nabla_X({\bf s}))'$. For
clarity we suppress the index $({\bf s})$:
\begin{align*}
f(\theta_X, \theta_Y) & = \int_0^\infty\int_0^\infty g(r_y\cos\theta_Y,
r_y\sin\theta_Y, r_x\cos\theta_X, r_x\sin\theta_X)r_yr_x dr_ydr_x\\
& = \int_0^\infty
r_x\frac{1}{\sqrt{(2\pi)^3|\Sigma|}}\exp(-\frac{1}{2}(ac-\tilde{A}^2)r_x^2)\\&\times\left(\int_0^\infty r_y
\frac{\sqrt{a}}{\sqrt{2\pi}}\exp(-\frac{a}{2}(r_y-\beta
r_x\cos\theta)^2)dr_y\right)dr_x\\
& = \frac{1}{\sqrt{(2\pi)^2|\Sigma|}}\int_0^\infty
\frac{1}{\sqrt{2\pi}}\left(\frac{\Phi(\tilde{A}r_x)\tilde{A}r_x/\sqrt{a}}{\sqrt{a}}+\frac{\phi(\tilde{A}r_x)}{a}\right)r_x\exp(-\frac{1}{2}(ac-\tilde{A}^2)r_x^2)dr_x\\
& = \frac{1}{\sqrt{(2\pi)^2|\Sigma|}}\Big(
  \frac{1}{\sqrt{2\pi a}} \frac{\tilde{A}}{\sqrt{a}}\int_0^\infty r_x\Phi(\tilde{A}r_x)r_x\exp(-\frac{1}{2}(ac-\tilde{A}^2)r_x^2) dr_x \\&+
\int_0^\infty
r_x\frac{1}{a\sqrt{2\pi}}\phi(\tilde{A}r_x)\exp(-\frac{1}{2}(ac-\tilde{A}^2)r_x^2)dr_x\Big)\\
& = \begin{cases}
\frac{1}{\sqrt{(2\pi)^2|\Sigma|}}\left(\frac{\tilde{A}}{a\sqrt{2\pi}}
  \left(\frac{\phi(0)\tilde{A}}{ac(ac-\tilde{A}^2)}+
  \frac{\sqrt{2\pi}}{(ac-\tilde{A}^2)^{3/2}}L(0,0,\sqrt{\frac{\tilde{A}^2}{ac}})\right)+\frac{\phi(0)}{a^2c\sqrt{2\pi}}
\right), & \text{$\tilde{A}>0$}\\
 \frac{1}{\sqrt{(2\pi)^2|\Sigma|}}\left(\frac{\tilde{A}}{a\sqrt{2\pi}}
  \left(\frac{\phi(0)\tilde{A}}{ac(ac-\tilde{A}^2)}+
  \frac{\sqrt{2\pi}}{(ac-\tilde{A}^2)^{3/2}}(0.5-L(0,0,\sqrt{\frac{\tilde{A}^2}{ac}}))\right)+\frac{\phi(0)}{a^2c\sqrt{2\pi}}
\right),&  \text{$\tilde{A}<0$}
\end{cases}
\end{align*}
where
$|\Sigma|=(\sigma_x^2\phi_x^2)^2(\sigma_y^2\phi_y^2)^2$,
$a=1/(\sigma_y^2\phi_y^2)$,
$c=(\sigma_y^2\phi_y^2+\beta^2\phi_x^2\sigma_x^2)/(\sigma_x^2\phi_x^2)$,
$\tilde{A}=\sqrt{a}\beta \cos(\theta_X-\theta_Y)$, and $L(0, 0, \rho)$
is the zero mean bivariate normal cdf with correlation $\rho$ and standard
deviations equal to 1 evaluated at $(0,0)'$.

Several integration steps were treated as the expected value of a
truncated normal distribution. An additional integration step required
integration by three parts, then computation of the probability that
the sum of a normal and a truncated normal are less than some
constant. This probability was available in closed form in the query
by \cite{Nelson1964} and is written in terms of the bivariate
normal cdf. Somewhat simplified, this gives the density provided in
the text:
\begin{align*}
f(\theta_X,\theta_Y) &= \begin{cases}
\frac{1}{a(2\pi)^{3/2}\sqrt{|\Sigma|}}\left(\frac{\tilde{A}^2\phi(0)}
{ac(ac-\tilde{A}^2)} + \frac{\sqrt{2\pi}\tilde{A}}{(ac-\tilde{A}^2)^{3/2}}
 L(0,0,\sqrt{\frac{\tilde{A}^2}{ac}}) +
 \frac{\phi(0)}{ac}\right),&\text{$\tilde{A}>0$}\\
\frac{1}{a(2\pi)^{3/2}\sqrt{|\Sigma|}}\left(\frac{\tilde{A}^2\phi(0)}
{ac(ac-\tilde{A}^2)} + \frac{\sqrt{2\pi}\tilde{A}}{(ac-\tilde{A}^2)^{3/2}}
 \left(0.5 - L(0,0,\sqrt{\frac{\tilde{A}^2}{ac}})\right) +
 \frac{\phi(0)}{ac}\right),&\text{$\tilde{A}<0$}
\end{cases}
\end{align*}


\section*{Acknowledgements}
  This research was supported in part by NSF DEB 0842465 and NSF
  CDI 0940671. The authors thank Thomas J. Leininger for
  computational assistance and James S. Clark for providing data and
  giving insight into the underlying ecological processes.

\bibliographystyle{elsarticle-harv}
\bibliography{References}{}

\end{document}